\documentclass[12pt]{amsart}
\allowdisplaybreaks
\usepackage{amssymb}
\usepackage{graphicx}
\usepackage{fullpage}
\usepackage[pdfencoding=auto,colorlinks,linkcolor=,citecolor=]{hyperref}
\usepackage{cite}
\usepackage{dsfont}
\usepackage{eucal}
\usepackage[all,cmtip]{xy}

\newcommand{\ol}{\overline}

\newcommand{\da}{\downarrow}

\newcommand{\Dim}{\mathsf{Dim}}

\newcommand{\ep}{\varepsilon}
\newcommand{\Ext}{\mathsf{Ext}}
\newcommand{\tExt}{\mathsf{E\widehat{\phantom{\dot{}}x}t}}

\newcommand{\Hom}{\mathsf{Hom}}
\newcommand{\HH}{\mathsf{H}}
\newcommand{\hH}{{\widehat \HH}}

\newcommand{\Ker}{\mathsf{Ker}}
\newcommand{\kD}{{\kk D}}
\newcommand{\kG}{{\kk G}}
\newcommand{\kH}{{\kk H}}
\newcommand{\kk}{\mathsf{k}}

\newcommand{\mmod}{\mathsf{mod}}

\newcommand{\Rad}{\mathsf{Rad}}

\newcommand{\Soc}{\mathsf{Soc}}

\newcommand{\Mod}{\mathsf{Mod}}
\newcommand{\stmod}{\mathsf{stmod}}
\newcommand{\StMod}{\mathsf{StMod}}

\newcommand{\ua}{\uparrow}

\newcommand{\bF}{\mathbb F}

\newcommand{\bZ}{\mathbb Z}

\newcommand{\cF}{\mathcal F}
\newcommand{\cG}{\mathcal G}

\newcommand{\cN}{\mathcal N}
\newcommand{\cNG}{\cN(\kk G)}

\newcommand{\hM}{\overline{M}}

\newcommand{\hU}{\widehat{U}}

\newcommand{\hV}{\widehat{V}}

\newcommand{\FG}{\cF\cG}
\newcommand{\FGG}{\FG(\kG)}

\renewcommand{\le}{\leqslant}
\renewcommand{\ge}{\geqslant}

\makeindex
\numberwithin{equation}{section}  

\theoremstyle{plain}
\newtheorem{lemma}[equation]{Lemma}
\newtheorem{theorem}[equation]{Theorem}
\newtheorem{proposition}[equation]{Proposition}
\newtheorem{corollary}[equation]{Corollary}

\theoremstyle{definition}

\theoremstyle{remark}

\author{David J. Benson} 
\address{Institute of Mathematics \\ 
University of Aberdeen \\ 
Aberdeen AB24 3UE \\ 
United Kingdom}
\email{d.j.benson@abdn.ac.uk}

\author{Jon F. Carlson}
\address{Department of Mathematics, University of Georgia, Athens GA
  30602, USA}

\subjclass[2010]{20C20, 20J06}

\title{Modules with finitely generated cohomology}

\begin{document}

\begin{abstract}
Let $G$ be a finite group and $\kk$ a field of characteristic $p$. It
is conjectured in a paper of the first author and John Greenlees that
the thick subcategory of the stable module category $\StMod(\kG)$
consisting of modules whose cohomology is finitely generated over
$\HH^*(G,\kk)$ is generated by finite dimensional modules and modules
with no cohomology.
If the centraliser of every element of order $p$
in $G$ is $p$-nilpotent, this statement follows from previous work.
Our purpose here is to prove this conjecture in two cases with non
$p$-nilpotent centralisers. The groups involved are
$\bZ/3^r\times\Sigma_3$ ($r\ge 1$) in characteristic three and $\bZ/2\times A_4$ in
characteristic two.

As a consequence, in these cases the bounded
derived category of $C^*BG$ (cochains on $BG$ with coefficients in
$\kk$) is generated by $C^*BS$, where $S$ is a Sylow $p$-subgroup of $G$.
\end{abstract}

\maketitle


\section{Introduction}

Let $G$ be a finite group and $\kk$ a field of characteristic $p$. Let
$\StMod(\kG)$ be the stable module category of $\kG$-modules and 
let $\stmod(\kG)$ be the subcategory of finitely generated
$\kG$-modules. These are triangulated categories, and the latter
consists of the compact objects in the former. 

Let $C^*BG$
the cochains on the classifying space of $G$ with coefficients in $\kk$.
In a recent paper~\cite{Benson/Greenlees:bg11}, it was proved that the
following two statements are equivalent:
\begin{enumerate}
\item[\rm (i)] The thick subcategory $\FGG$ of $\StMod(\kG)$ consisting of
  modules $M$ such that the non-negative Tate cohomology
  $\hH^{\ge 0}(G,M)$ is finitely generated as a module
  over $\HH^*(G,\kk)$ is equal to the thick subcategory generated by the
  finite dimensional modules and the modules $M'$ satisfying 
$\hH^{\ge 0}(G,M')=0$.
\item[\rm (ii)] The bounded derived category of $C^*BG$ is generated
  by $C^*BS$, where $S$ is a Sylow $p$-subgroup of $G$.
\end{enumerate}
It was also proved that if the
centraliser in $G$ of every element of order $p$ is $p$-nilpotent, then
statements (i) and (ii) are true. However, in this case the reason why
statement (i) holds is that, if $\hH^{\ge 0}(G,M)$ is
finitely generated over $\HH^*(G,\kk)$, then $M$ decomposes as a direct sum
$M=M'\oplus M''$ such that $\hH^*(G,M')=0$ and $M''$ is finite
dimensional. This follows from the main theorems
of~\cite{Benson:1995a,Benson/Carlson/Robinson:1990a}, but
is false more generally, as is easy to see in examples.

Let $\cNG$ denote the thick subcategory of no-cohomology modules, meaning
the $\kG$-modules $M$ such that $\hH^*(G,M)=0$. Statement (i) says that 
$\FGG$ is generated as a triangulated subcategory 
by objects in the two subcategories $\cNG$ and 
$\stmod(\kG)$ of $\StMod(\kG)$. 

The purpose of this paper is to prove statement (i) in two of the 
simplest cases in which the centralizer of some
$p$-element of $G$ is not $p$-nilpotent. The cases are 
\begin{itemize}
\item[(A).] $G \cong \bZ/3^r \times \Sigma_3$ for some $r \ge 1$, 
in characteristic $p=3$, and
\item[(B).] $G \cong \bZ/2 \times A_4$ in characteristic $p=2$.
\end{itemize}
 More precisely, we prove that the following 
refinement of statement (i) above holds.

\begin{theorem}\label{thm:tech}
Let $G$, $r$ and $p$ be as in case $(A)$ or $(B)$, as above. Let 
$\kk$ a field of characteristic $p$.
Given a $\kG$-module $M$ such that $\HH^i(G,M)$ is finite dimensional
for each $i \ge 0$,
there are submodules $M_1\le M_2 \le M$ such that 
$\hH^*(G,M_2/M_1)=0$, and $M_1$ and $M/M_2$
are finite dimensional. 
In particular,
$\hH^{\ge 0}(G,M)$ is finitely generated as a module over $\hH^*(G,\kk)$.
\end{theorem}

This leads directly to our main theorem. 

\begin{theorem}\label{th:main}
Suppose that $G$ and $\kk$ satisfy one of the following hypotheses:
\begin{enumerate}
\item
The centraliser of every element of order $p$ is
$p$-nilpotent,
\item $G=\bZ/3^r\times \Sigma_3$, $r \ge 1$ and $\kk$ has
characteristic three, or
\item $G=\bZ/2\times A_4$ and $\kk$ has
characteristic two.
\end{enumerate}
Then $\FGG$ is generated by $\cNG$ and $\stmod(\kG)$, as a thick
subcategory of $\StMod(\kG)$.
\end{theorem}

The case where the centralisers are $p$-nilpotent is an easy
consequence of the main theorems
of~\cite{Benson:1995a,Benson/Carlson/Robinson:1990a}.
The other two cases follow directly from the previous theorem.

After some preliminaries, the case $G=\bZ/3^r\times\Sigma_3$ is dealt
with in Section~\ref{se:Z/3xS3}, and the case $G=\bZ/2 \times A_4$ is
proved in Section~\ref{se:Z/2xA4}.\bigskip

\noindent
{\bf Acknowledgements.} Both authors were supported by the Hausdorff
Institute in Bonn while part of this research was carried out. They also thank
John Greenlees and Henning Krause for numerous conversations about
this work.


\section{Some preliminaries}

In this section, we recall some required facts concerning 
group cohomology and the structure of no-cohomology modules.
Throughout the paper, $G$ is a finite group and $\kk$ denotes a 
field of characteristic $p > 0$. 
Let $\Mod(kG)$ denote the category of all $\kG$-modules. The 
stable category $\StMod(\kG)$ has the same objects as $\Mod(\kG)$,
but the morphisms between two objects is the quotient of the group 
of $\kG$-module homomorphisms by the subgroup of those that factor through 
a projective module. Let $\mmod(\kG)$ be the subcategory of $\Mod(\kG)$
consisting of finitely generated modules. Then $\stmod(\kG)$ is the 
subcategory of $\StMod(\kG)$ of objects isomorphic to finitely 
generated modules. The categories $\stmod(\kG)$ and $\StMod(\kG)$ 
are tensor triangulated categories with the tensor product given
by the Hopf algebra structure on $\kG$. The triangulated structure is described in detail in
the first chapter of Happel~\cite{Happel:1988a}, and we shall not be
using the tensor structure.

As in the introduction, $\FGG$ and $\cNG$ are the thick subcategories of 
$\StMod(\kG)$ consisting respectively of the modules with finitely
generated cohomology and the no-cohomology modules. 

\begin{lemma}\label{le:BCRo}
There exists a constant $n_p(G)$, depending only on $G$ and $p$, such
that the following hold for any $\kG$-module $M$.
\begin{enumerate}
\item If $\hH^i(G,M)$ vanishes for $n_p(G)$ successive values of $i$
then $M$ is in $\cNG$, and $\hH^i(G,M) = \{0\}$ for all $i\in\bZ$.
\item If $\hH^i(G,M)$ is finite dimensional for $n_p(G)$ successive
values of $i$ then $M$ is in $\FGG$ and $\hH^i(G,M)$ is finite 
dimensional for all $i\in\bZ$.
\end{enumerate}
\end{lemma}
\begin{proof}
Statement (1) is proved in~\cite{Benson/Carlson/Robinson:1990a}, and
statement (2) is proved in~\cite{Benson/Greenlees:bg11} using the same
technique.
\end{proof}

We recall that if $M$ is a module with finitely generated
cohomology that has no projective submodules, then $\Hom_{\kG}(\kk,M)= 
\HH^0(G,M) = \hH^0(G,M)$ has finite dimension. 

The following is standard Tate duality. 

\begin{proposition} \label{pr:Tate-duality}
If $N$ is a finitely generated $\kG$-module and $M$ is an arbitrary
$\kG$-module, then we have
\[ \tExt_\kG^{n}(M,N) \cong \Hom_\kk(\tExt_\kG^{-n-1}(N,M),\kk). \]
\end{proposition}

\begin{corollary}\label{co:Hom(M,k)}
If $\HH^i(G,M)$ is finite dimensional for $i\ge 0$, then the spaces
$\Hom_\kG(\kk,M)$ and
$\Hom_{\kG}(M,\kk)$ are also finite dimensional.
\end{corollary}
\begin{proof}

  If $\HH^i(G,M)$ is finite dimensional for $i\ge 0$ then 
  $\HH^0(G,M) \cong \Hom_{\kG}(\kk,M)$ is finite dimensional. This
  implies that $M$
  has only a finite number of copies of the projective cover $P_\kk$
  of the trivial module $\kk$ as summands (since $P_\kk$ is also
  the injective hull of $\kk$). Also 
$\hH^i(G,M)$ is finite dimensional for $i\ge 0$. Applying
  Lemma~\ref{le:BCRo}\,(2), it follows that $\hH^i(G,M)$ 
is finite dimensional for all
  $i\in\bZ$. In particular, $\hH^{-1}(G,M)$ is finite dimensional.
  Applying Proposition~\ref{pr:Tate-duality} with $N=\kk$ and $n=0$,
  we deduce that $\tExt^0_\kG(M,\kk)$ is finite dimensional. Since $M$
  has only a finite number of copies of $P_\kk$ as summands, this
  implies that $\Hom_\kG(M,\kk)$ is finite dimensional.
\end{proof}

To make use of this, we have the following.

\begin{proposition}\label{le:O^pG}\ 
Let $M$ be a $\kG$-module. 
\begin{enumerate} 
\item[\rm (1)]  If $\Hom_\kG(\kk,M)$ is
finite dimensional then $M$ has a finite dimensional submodule $M'$
such that $\Hom_\kG(\kk,M/M')=0$.
\item[\rm (2)] If $\Hom_\kG(M,\kk)$ is finite
dimensional, then $M$ has a submodule $M''$ such that $M/M''$ is finite
dimensional and $\Hom_\kG(M'',\kk)=0$.
\end{enumerate}
Moreover, $M'$ and $M/M''$ can be chosen so that all of their 
composition factors are isomorphic to $\kk$.
\end{proposition}

\begin{proof}
Let $I_M$ be an injective hull of $M$. Then
$\Hom_\kG(\kk,I_M)\cong\Hom_\kG(\kk,M)$ is finite dimensional, and so
$I_M$ has a finite number of copies of the injective hull $I_\kk$ of $\kk$ as a
summand. Let $I_M=U\oplus I$ where $U$ is a direct sum of a finite
number of copies of $I_\kk$ and $\Hom_\kG(\kk,I)=0$. Then every submodule
of $I_\kk$ whose composition factors are all copies of $\kk$ is contained
in $U$. It follows that every submodule of $M$ whose composition
factors are all copies of $\kk$ is finite dimensional. There is a
unique maximal one, and we can take this as $M'$.

Dually, using
projective covers, there is
a unique minimal submodule of $M$ such that all the composition
factors of the quotient are isomorphic to $\kk$. We can take this as
$M''$, and then the quotient $M/M''$ is finite dimensional.
\end{proof}

A combination of the two results above yields the following.

\begin{corollary} \label{cor:core}
Suppose that $M$ is a $kG$-module having finitely generated cohomology and
no projective submodules. Then $M$ has submodules $M_1 \le M_2 \le M$
such that
\begin{enumerate}
\item the modules $M_1$ and $M/M_2$ have finite dimension with all of their
composition factors isomorphic to $\kk$, and
\item $\Hom_\kG(M_2/M_1,\kk) = \{0 \} = \Hom_\kG(\kk, M_2/M_1)$.
\end{enumerate}
\end{corollary}

\begin{proof} 
Use part (1) of the Proposition \ref{le:O^pG} to find $M_1 = M'$ with
$\Hom_\kG(\kk, M/M_1) = \{0 \}$. Then, by part (2), there is a 
submodule $M_2/M_1$ of $M/M_1$ such that $\Hom_\kG(M_2/M_1, \kk) = \{0 \}$.
Now note that $\Hom_\kG(\kk, M_2/M_1)$ is a subsubspace of 
$\Hom_\kG(\kk, M/M_1)$, and hence it is also zero.
\end{proof}

We also need some more specific facts pertinent to the representation
theory of the groups we are analyzing in the next two sections. 

\begin{lemma}\label{le:Hom-k}
Let $H$ be a normal subgroup of $p$-power index in a finite group $G$
and let $\kk$ be a field of characteristic $p$. If $M$ is a
$\kG$-module such that
$\Hom_\kG(\kk,M)=0$, then $\Hom_\kH(\kk,M{\da_H})=0$. In the other 
direction, if $\Hom_\kG(M,\kk)=0$, then $\Hom_\kH(M{\da_H},\kk)=0$.
\end{lemma}
\begin{proof}
Suppose first that $\Hom_\kG(\kk,M)=0$. We claim that if $N$ is any
module having all composition factors isomorphic to $\kk$ then
$\Hom_\kG(N,M)=0$. The proof is by induction on the number of
composition factors. The case of just one composition factor 
is by hypothesis, and the
inductive step follows from the left exactness of $\Hom$.
Now take $N=\kk_H{\ua^G}$. This only has $\kk$ as composition factors
because $G/H$ is a finite $p$-group. Thus by what we just proved, we have
$0= \Hom_\kG(\kk{\ua^G},M) \cong \Hom_\kH(\kk,M{\da_H})$, 
proving this case of the lemma.
The case where $\Hom_\kG(M,\kk)=0$ is dual.
\end{proof}

Our examples both have a normal Sylow $p$-subgroup $P$ with cyclic
quotient $\langle t\rangle$. So we shall use generators of the radical
of $\kk P$ consisting of eigenvectors of conjugation by $t$, and bases of
$\kG$-modules also consisting of
eigenvectors of $t$. We make use of the following rather obvious
general principle.

\begin{proposition}\label{pr:eigenvalues}
Suppose that $\kk$ is a field of characteristic $p$, and the group
$G$ is a semi-direct product $P\rtimes\langle t\rangle$ where $P$
is a Sylow $p$-subgroup and $t$ has order coprime to $p$.
Let $M$ be a $\kG$-module. If $m\in M$ is an eigenvector of $t$ with
eigenvalue $\zeta\in\kk^\times$ and $X$ is an element of $\kk P$ with
$tX=\eta Xt$ with $\eta\in\kk^\times$, then $Xm$ is an eigenvector of
$t$ with eigenvalue $\eta\zeta$.
If $m=Xm'$
for some $m'\in M$, then there exists an eigenvector $m_1$ of $t$
with eigenvalue $\eta^{-1}\zeta$ such that $m = xm_1$.
\end{proposition}

\begin{proof}
We have $t(Xm)=(\eta Xt)m=(\eta X)(\zeta m)=(\eta\zeta)(Xm)$, proving
the first part. For the second part, write $m'$ as a sum of
eigenvectors of $t$, say $m'=m_1+\dots+m_q$, where $m_1$ has
eigenvalue $\eta^{-1}\zeta$ and the rest have other eigenvalues. Then
$m=Xm_1+\dots + Xm_n$ is a sum of eigenvectors. The eigenvalue of $t$
on $Xm_1$ is $\eta(\eta^{-1}\zeta)=\zeta$, and the rest have other
eigenvalues. So $m=Xm_1$, and $Xm_i=0$ for $i\ge 2$.
\end{proof}

We use the following transitivity lemma on several occasions. 

\begin{lemma}  \label{lem:trans}
Suppose that $M$ is a $\kG$-module and that $L = M_2/M_1$ for submodules
$M_1 \le M_2 \le M$. Suppose that $N = L_2/L_1$ for  submodules 
$L_1 \le L_2 \le L$. Then there exist submodules $M_1' \le M_2' \le M$
such that $N \cong M_2'/M_1'$. 
\end{lemma} 

\begin{proof}
Since $L_1 \le M_2/M_1$, there is a submodule $M_1' \le M_2$ such that 
$M_1'/M_1 = L_1$. Likewise there is $M_2'$ such that $L_2 = M_2'/M_1$. 
A standard isomorphism theorem says that $N \cong M_2'/M_1'$. 
\end{proof}

\begin{proposition}\label{pr:res-nocoho}
Suppose that $H$ is a subgroup with $p$-power index in a finite group $G$.
Let $M$ be a $\kG$-module with the 
properties that the restriction 
$M{\da_H}$ is projective and that $\Hom_{\kH}(\kk,M{\da_H}) = \{0\}$. 
Then $M$ is a no-cohomology module.
\end{proposition}

\begin{proof}
Let $(P_*, \ep)$ be an injective $\kk(G/H)$-resolution of $\kk$. Note 
that because $G/H$ is a $p$-group, it is a resolution of free modules,
meaning the every $P_i$ is a direct sum of copies of 
$\kk (G/H)$. By inflation we may regard $\kk (G/H)$ as a 
$\kG$-module on which $H$ acts trivially.  We may also do the 
same for the resolution $P_*$. Then
tensoring  with $M$, we get the sequence 
\[
\xymatrix{
0 \ar[r] & M \ar[r]^{1 \otimes \ep \quad} & M \otimes P_0 \ar[r] &
M \otimes P_1 \ar[r] & M \otimes P_2 \ar[r] & \cdots
}
\]
As a $\kG$-module $\kk (G/H) \cong \kk{\da_H}{\ua^G}$.
By Frobenius reciprocity, 
\[
M \otimes \kk{\da_H}{\ua^G} \cong M{\da_H}{\ua^G}.
\]
Because $M{\da_H}$ is projective, we have that
$M \otimes \kk(G/H)$ is projective and hence 
$M\otimes P_i$ is projective as a $\kG$-module for all $i$.
Thus, the above sequence is an injective resolution of $M$.

Moreover, 
\[
\Hom_{\kG}(\kk, M{\da_H}{\ua^G}) \cong 
\Hom_{\kH}(\kk, M{\da_H}) = \{0\}.
\]
Consequently, $\Hom_{\kG}(\kk, M \otimes P_i) = \{0\}$ and 
$\HH^i(G, M) = \{0\}$ for all $i\ge 0$. Finally, applying Lemma~\ref{le:BCRo},
it follows that $\hH^i(G,M)=\{0\}$ for all $i\in\bZ$.
\end{proof}


\section{\texorpdfstring{The case $G=\bZ/3^r \times \Sigma_3$}
  {The case G = ℤ/3ʳ × Σ₃}}\label{se:Z/3xS3}

Throughout this section let $\kk$ be a field of 
characteristic three and let $r \ge 1$ be an integer. 
Let $G=C\times D$ where $C=\langle g \mid g^{3^r} =1\rangle\cong \bZ/3^r$ and
\[ D=\langle h,t\mid h^3=1, t^2=1, tht=h^{-1}\rangle \cong \Sigma_3. \]
Set $Z=g-g^2$ and
$Y=h-h^2$. Then
\[ \kG=\kk\langle Z,Y,t\mid Z^{3^r}=Y^3=0,\ ZY=YZ,\ tZ=Zt,\ tY=-Yt\rangle. \]
Up to isomorphism, there are two simple $\kG$-modules, both having dimension one, 
which we write as $\kk$ and $\ep$, with $t$ acting as $1$ on $\kk$ and
as $-1$ on $\ep$.

In the case that $r=1$, the structure of the 
projective indecomposable $\kG$-modules may be
pictured as follows.
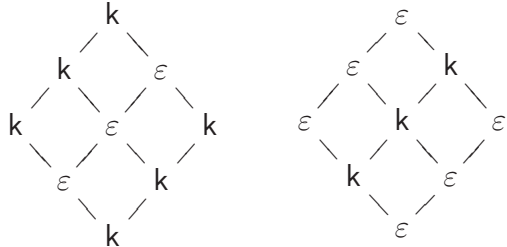
\[ 
\xymatrix@=2mm{
&&\kk\ar@{-}[dl]\ar@{-}[dr]\\
&\kk\ar@{-}[dl]\ar@{-}[dr]&&\ep\ar@{-}[dl]\ar@{-}[dr]\\
\kk\ar@{-}[dr]&&\ep\ar@{-}[dl]\ar@{-}[dr]&&\kk\ar@{-}[dl]\\
&\ep\ar@{-}[dr]&&\kk\ar@{-}[dl]\\
&&\kk}\qquad
\xymatrix@=2mm{
&&\ep\ar@{-}[dl]\ar@{-}[dr]\\
&\ep\ar@{-}[dl]\ar@{-}[dr]&&\kk\ar@{-}[dl]\ar@{-}[dr]\\
\ep\ar@{-}[dr]&&\kk\ar@{-}[dl]\ar@{-}[dr]&&\ep\ar@{-}[dl]\\
&\kk\ar@{-}[dr]&&\ep\ar@{-}[dl]\\
&&\ep
} 
\]
In the diagrams, a node represents a basis element on which $t$ acts by
multiplication by either $1$ or $-1$ as indicated. A line going down left 
indicates that the lower left basis element is $Z$ times the basis element
at the upper right end of the line. Going down right is multiplication by $Y$.

If $r > 1$, then the diagrams must be extended down left to $3^r$ rows, of
$[\kk, \ep, \kk]$ or $[\ep, \kk, \ep]$. 
In any case, the $\Ext^1$ groups between simple modules are all one
dimensional.\medskip

We can explicitly list the indecomposable $\kD$-modules. The simple
modules are $\kk$ and $\ep$. The
projective covers $P_\kk$ of $\kk$ and $P_\ep$ of $\ep$ are  
uniserial $[\kk,\ep,\kk]$ and $[\ep,\kk,\ep]$, repsectively. 
All indecomposable modules are uniserial. The non-projective indecompoable 
modules are the simple modules together with the length two
uniserial modules $[\kk,\ep]$ and $[\ep,\kk]$. So the list is:
\begin{equation}\label{eq:D}
\kk\qquad \ep\qquad
\vcenter{\xymatrix@=2mm{\kk\ar@{-}[d]\\\ep}}\qquad
\vcenter{\xymatrix@=2mm{\ep\ar@{-}[d]\\\kk}}\qquad \quad 
P_{\kk} = 
\vcenter{\xymatrix@=2mm{\kk\ar@{-}[d]\\\ep \ar@{-}[d]\\\kk}}\qquad
\quad P_{\ep} =
\vcenter{\xymatrix@=2mm{\ep\ar@{-}[d]\\\kk \ar@{-}[d]\\\ep}}\ .
\end{equation}

The purpose of this section is to prove Theorem \ref{thm:tech} in the 
case that group $G$ is given as above (case (A) in the Introduction). 
The proof proceeds by a series of reductions, the last of which is a
mathematical induction argument. A basic principle is the transitivity
statement expressed in Lemma \ref{lem:trans}. In particular, supppose
that $M$ is a $\kG$-module  having submodules $L_1 \le L_2 \le M$
such that both $L_1$ and $M/L_2$ are finite dimensional. Then 
if $L_2/L_1$ satisfies the conditions of Theorem \ref{thm:tech},
then so also does $M$. 

\begin{proof}[Proof of Theorem \ref{thm:tech} in case (A)]
Suppose that $M$ is a $\kG$-module such that $\HH^i(G,M)$ is finite
dimensional for all $i$. Our objective is to show that 
$M$ satisfies the conditions stated in Theorem~\ref{thm:tech}. 
By Corollary~\ref{cor:core} and the transitivity  we may 
(and do) assume that
\[ \Hom_{\kG}(M,\kk) = \{0\} = \Hom_{\kG}(\kk,M). \]

{\bf Step 1:} We have that $\Hom_{\kD}(\kk, M{\da_D})=0= 
\Hom_{\kD}(M{\da_D},\kk)$. 
The restricted module $M{\da_D}$ is a direct sum of copies of $\ep$
and the projective cover $P_\ep$, which is the uniserial 
module $[\ep,\kk,\ep]$. The number of copies of $\ep$ in a 
direct sum decomposition of $M{\da_D}$ is finite.\medskip

The first statement follows from Lemma~\ref{le:Hom-k},
since $D$ has index $3^r$ in $G$.
From the classification described above, the only
indecomposable $\kD$-modules with no homomorphisms to or from $\kk$
are $\ep$ and the projective cover $P_\ep$.
Thus the second statement follows.

If the number of summands isomorphic to $\ep$ in $M{\da_D}$ 
were infinite, then $\Ext^1_{\kk D}(\kk,M{\da_D})$ 
would be infinite dimensional.
This is isomorphic to $\Ext^1_{\kG}(\kk_D{\ua^G},M)$. The module
$\kk_D{\ua^G}$ is uniserial of length $3^r$ with all composition 
factors isomorphic to $\kk$. Hence, by the long 
exact sequence in $\Ext$, the degree one cohomology $\HH^1(G,M)\cong
\Ext^1_{\kG}(\kk,M)$ would be infinite dimensional, contradicting our
assumptions on $M$.\medskip

{\bf Step 2:}  Let $d(M)$ denote the number of copies of the 
simple module $\ep$ that occur in a direct sum decomposition of
$M{\da_D}$. This number is well defined
since it is the dimension of the $\kk D$-module $\Soc(M{\da_D})/Y^2M{\da_D}$
which is finite dimensional.\medskip

The remainder of the proof of Theorem \ref{thm:tech} is by 
induction on $d(M)$. The next step is to prove the base case.\medskip

{\bf Step 3:} If $d(M) = 0$, then $M$ is a no-cohomology module.\medskip

The restriction $M{\da_D}$ is a projective $\kk D$-module satisfying
$\Hom_{\kk D}(\kk,M{\da_D})=0$, and
so applying Proposition~\ref{pr:res-nocoho}, 
it follows that $M$ is a no-cohomology module.\medskip

{\bf Step 4:}  
Let $L=L(M) = \Ker(Y^2,M) = \{m \in M \mid Y^2m = 0\}$. We may assume that 
\begin{enumerate}
\item $\Soc(M) \le Y^2M\le \Soc_{\kD}(M{\da_D})$, and that
\item $\Rad_{\kD}(M{\da_D}) \le L(M) \le \Rad_{kG}(M)$.
\end{enumerate}
Proof of (1): For the first containment,
if $m \in \Soc(M)$ is not in $Y^2M$, then $m$ generates a $\kG$-submodule 
$M_1 \cong \ep$ that is not contained in any projective $kD$-submodule. 
This implies that $d(M/M_1) < d(M)$, and we are done by induction.
The second containment follows from the fact that the restriction of $M$
to $\kD$ is a direct sum of copies of $\ep$ and $P_\ep$.

The proof of (2) is dual 
to that of statement (1).\medskip

{\bf Step 5:} 
Next choose an element $m \in L = L(M)$ having the following properties. 
\begin{itemize}
\item[(i)] $tm = -m,$
\item[(ii)] $m \not\in Y^2M,$
\item[(iii)] $Zm \in Y^2M$, and 
\item[(iv)] $\Dim(\kG m)$ is minimal among such choices. 
\end{itemize}
Note that there must be such an element because $d(M) \neq 0$. That is, 
we can choose $m$ to be a generator of some direct summand of $M{\da_D}$
that is isomorphic to $\ep$. Then the first two conditions are satisfied. 
It is clear that a choice can be made that satisfies the remaining two
conditions.  

Next find an element $u \in M$ such that $Y^2u = Zm$. 
We claim that $u \not\in \Rad_{kG}(M).$  Otherwise, there exist $v \in M$
such that $Zv \equiv u$ modulo $YM$. Then $ZY^2v=Y^2u=Zm$, so
$Z(m-Y^2v) = 0$.
Also $Y(m-Y^2v)=0$, so $m-Y^2v\in\Soc_{\kG}(M)\le Y^2M$ by statement
(1) of Step 4. Hence $m\in Y^2M$, contradicting condition (ii).\medskip

{\bf Step 6:}
Let $U$ denote the $\kG$-submodule of $\hM = M/L$ generated by 
$\ol{u}=u+L$. We claim that there is a submodule $V$ such that 
$\hM \cong U \oplus V$ as $kG$-modules. The proof goes as follows. 

Note first that $Y$ annihilates $\hM$, and $t$ acts by 
multiplication by $-1$ on the entire module. So $\hM$ may 
be regarded as a module over the algebra $A = \kk[Z]/(Z^{3^r}) 
= \kk\langle g \rangle$.
As such, it is a direct sum $\hM = \bigoplus_{i \in I} M_i$ for uniserial 
modules $M_i$ and some
indexing set $I$. Each $M_i$ is generated by an element $b_i$ and 
has a $\kk$-basis $b_i, Zb_i, \dots, Z^{n_i-1}b_i.$ That is, 
$n_i$ is the length of $M_i$. 

Let $n$ be the length of $U$, so that $U$ has a basis $\ol{u}$,
$Z\ol{u}, \dots, Z^{n-1}\ol{u}$. It follows that the collection
$m, Zm, \dots, Z^nm$ is a basis for $\kG m$. 

Now write $\ol{u} = \sum_{i \in I} \alpha_ib_i$ for $\alpha_i \in A$, where
all but a finite number of $\alpha_i$'s are zero. We assert that 
there is some $i \in I$ such that $\alpha_i \not\in \Rad(A)$ and $n_i = n$.
For suppose not. Then for every $i$ with $\alpha_i \neq 0$, we must 
have that either $n_i < n$ or $\alpha_i \in \Rad(A)$, meaning that 
$\alpha_i = Z\alpha_i'$ for some $\alpha_i' \in A$. It follows that 
\[
\ol{u} = Z\ol{u}_1 +\ol{u}_2 \quad \text{ where } \quad Z^{n-1}\ol{u}_2 = 0 
\]

Let $u_1$ and $u_2$ be elements of $M$ such that $\ol{u}_1 = u_1+L$ and 
$\ol{u}_2 = u_2+L$. Consider the element $m'= m -Y^2u_1$. It is easy 
to check that $m'$ satisfies conditions (i), (ii) and (iii) of 
Step 5. But then 
\[
Zm' = Zm-ZY^2u_1= Y^2u-Y^2Zu_1 = Y^2u_2 \quad \text{ and } Z^nm' = Z^{n-1}Y^2u_2 = 0.
\]
However, this violates assumption (iv) on the minimality of the dimension
of the module generated by $m$. 

So we have that there exists some $i \in I$ such that 
$\alpha_i \not\in \Rad(A)$ and $n_i = n$.
Thus, $M_i \cong U$. Let $V$ be the submodule that is the direct sum of 
all $M_j$ for $j \neq i$. It is now easy to see that $\hM$ is the 
direct sum of $U$ and $V$, as claimed.\medskip

{\bf Step 7:}
Let $N = \Ker(Y,M) = \{m \in M \mid Ym = 0 \}\le L$. We claim that 
$M/N = \hU \oplus \hV$ where $\hU/Y\hU = U$ and $\hV/Y\hV = V$,
and multiplication by $Y$ induces isomorphisms
\[
\xymatrix{
\hU/Y\hU \ar[r] & Y\hU & \text{ and } & \hV/Y\hV \ar[r] & Y\hV
}
\]
To see this, we first notice that $(M/N){\da_D}$ is a direct sum of 
copies of the uniserial module $[\ep, \kk]$. The action of $t$ makes
this the direct sum of two eigenspaces
\begin{align*}
W_{-1} = (t-1)(M/N)&=\Ker(t+1,M/N), \\
W_{1} = (t+1)(M/N)&=\Ker(t-1,M/N)=Y(M/N) 
\end{align*}
corresponding to eigenvalues $-1$ and $1$
respectively. These are $\kk\langle g,t \rangle$-modules, and the
action of $Y$ maps the first isomorphically to the second.
In particular, As $\kk\langle g,t \rangle$-modules, 
$Y(M/N) = W_1 = L/N$ and $W_{-1} \cong M/L$.
Letting $\widetilde U$ and $\widetilde V$ be the inverse images of $U$ and $V$
under the map $M/N\to M/L$, we then set $\hU=(t-1)\widetilde U\oplus
Y\widetilde U$ and $\hV=(t-1)\widetilde V\oplus Y\widetilde V$ as internal direct sums of 
$\kk\langle g,t \rangle$-modules, with the obvious action of $Y$.\medskip

{\bf Final Step 8: } Let $M_1 = \kG m$ and let $M_2$ be the inverse image
of $\hV$ under the quotient map $M \to M/N$. 
Then $M/M_2 \cong \hU$ is finite dimensional, as is $M_1$.
We have that $M_1$ is spanned by $m, Zm = Y^2u, \dots, Z^nm = Z^{n-1}Y^2u$, 
while $M_2/N$ is a complement of the subspace spanned by the images of 
the elements of the set $\{Z^iu, Z^iYu | \ 0 \leq i \leq n-1\}$ 
in $M/N$, and $Y^2M_2$
is a vector space complement to the subspace spanned by 
$Y^2u, \dots,  Z^{n-1}Y^2u$ in $Y^2M$. 
Thus, we see from the restriction $M{\da_D}$, that 
$d(M_2/M_1) = d(M)-1$. Therefore, 
by induction on $d(M)$, we are done in this case. 
\end{proof}


\section{\texorpdfstring{The case $G=\bZ/2 \times A_4$}
  {The case G = ℤ/2 × A₄}}\label{se:Z/2xA4}

Let $\kk$ be a field of characteristic two 
containing $\bF_4=\{0,1,\omega,\bar\omega\}$, with $\omega^3 =1$. 
Let $G=C\times A$ where $C=\langle z\mid z^2=1\rangle$ and 
\[ 
A = \langle g,h,t\mid g^2=h^2=1,\ gh=hg,\ gt=th,\ 
ht=tgh,\ ght=tg,\ t^3=1\rangle. 
\]
Set $X= g+\omega h+\bar\omega gh$, $Y= g+\bar\omega h+\omega gh$, $Z=1+z$. Then
\[ \kG = \kk[X,Y,Z]/(X^2,Y^2,Z^2) \rtimes \langle t\mid t^3=1\rangle \]
with $tX=\omega Xt$, $tY=\bar\omega Yt$, $tZ=Zt$.

Up to isomorphism, $\kG$ has three simple modules, all one dimensional, 
which we shall denote $\kk$, $\omega$ and $\bar\omega$
on which $t$ acts by multiplication by the corresponding scalar.
If $m$ is a generator for one of these simple modules 
then $tm=m$, respectively $tm=\omega m$, $tm=\bar\omega m$.

The structure of the projective indecomposable $\kG$-modules may be
pictured as follows.
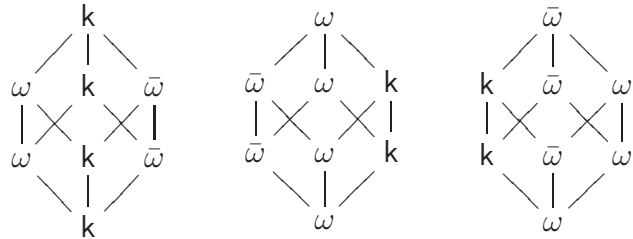
\[ \xymatrix@=4mm{&\kk\ar@{-}[dl]\ar@{-}[d]\ar@{-}[dr]\\
\omega\ar@{-}[d]\ar@{-}[dr]&\kk\ar@{-}[dl]\ar@{-}[dr]&
\bar\omega\ar@{-}[dl]\ar@{-}[d]\\
\omega\ar@{-}[dr]&\kk\ar@{-}[d]&\bar\omega\ar@{-}[dl]\\
&\kk}\qquad
\xymatrix@=4mm{&\omega\ar@{-}[dl]\ar@{-}[d]\ar@{-}[dr]\\
\bar\omega\ar@{-}[d]\ar@{-}[dr]&\omega\ar@{-}[dl]\ar@{-}[dr]&
\kk\ar@{-}[dl]\ar@{-}[d]\\
\bar\omega\ar@{-}[dr]&\omega\ar@{-}[d]&\kk\ar@{-}[dl]\\
&\omega}\qquad
\xymatrix@=4mm{&\bar\omega\ar@{-}[dl]\ar@{-}[d]\ar@{-}[dr]\\
\kk\ar@{-}[d]\ar@{-}[dr]&\bar\omega\ar@{-}[dl]\ar@{-}[dr]&
\omega\ar@{-}[dl]\ar@{-}[d]\\
\kk\ar@{-}[dr]&\bar\omega\ar@{-}[d]&\omega\ar@{-}[dl]\\
&\omega}
\]
The lines going down left, down right and down indicate multiplication 
by $X$, $Y$ and $Z$, respectively. 
Again, all the $\Ext^1$ groups between simple modules are one dimensional.

Here we begin the Proof of Theorem \ref{thm:tech} in case (B).
For this purpose we fix a $\kG$-module such that $\HH^i(G,M)$ is 
finite dimensional for all $i \ge 0$. The goal is
to show that we may remove finite dimensional pieces 
from the top and bottom of $M$,
leaving a no-cohomology module. We follow many of the steps of the 
last section, but with some significant differences. 
Assume the notation at the beginning of the section. 

\begin{lemma}\label{le:HomA}
Without loss of generality, we assume that $M$ has no projective summands,
and that $\Hom_{\kG}(\kk,M)= \{0\}$ and $\Hom_{\kG}(M,\kk)= \{0\}$.
As a consequence, $\Hom_{\kk A}(\kk,M{\da_A})=0$ and 
$\Hom_{\kk A}(M{\da_A},\kk)=0$.
\end{lemma}

\begin{proof}
The statement about projective summands is obvious. The next statement 
follows from \ref{cor:core} and the transitivity lemma \ref{lem:trans}.
The last statement follows from Lemma \ref{le:Hom-k}.
\end{proof}

\begin{lemma}\label{le:MdaA}
The restriction $M{\da_A}$ is a direct sum of copies of 
the simple modules
$\omega$, $\bar\omega$, their projective covers $P_\omega$ and
$P_{\bar\omega}$, and two length two uniserial modules
$[\omega,\bar\omega]$ and $[\bar\omega,\omega]$. 
\end{lemma}

\begin{proof}
The list of possible summands follows from Lemma~\ref{le:HomA} and the fact that
these are the indecomposable $\kk A$-modules $M$ satisfying
$\Hom_{\kk A}(\kk,M)=0$ and $\Hom_{\kk A}(M,\kk)=0$.
\end{proof}

We display the list of possible summands of $M{\da_A}$ described in
Lemma~\ref{le:MdaA} as follows: 
\begin{equation}\label{eq:A}
\omega\qquad\bar\omega\qquad
\vcenter{\xymatrix@=3mm{ &\omega\ar@{-}[dl]\\\bar\omega}}\qquad
\vcenter{\xymatrix@=3mm{\bar\omega\ar@{-}[dr]\\ &\omega}}\qquad
\vcenter{\xymatrix@=3mm{&\omega\ar@{-}[dl]\ar@{-}[dr]\\
    \bar\omega\ar@{-}[dr]&&\kk\ar@{-}[dl]\\&\omega}}\qquad
\vcenter{\xymatrix@=3mm{&\bar\omega\ar@{-}[dl]\ar@{-}[dr]\\
    \kk\ar@{-}[dr]&&\omega\ar@{-}[dl]\\&\bar\omega}}
\end{equation}

\begin{lemma} 
The number of non-projective summands in $M{\da_A}$ is finite.
\end{lemma}

\begin{proof}
The lemma is a consequence of the fact that the possible non-projective
summands all have cohomology:
\begin{align*}
\Ext^1_{\kG}(\kk_A{\ua^G},\omega)&\cong \HH^1(G,\omega)\cong \kk &
\Ext^1_{\kG}(\kk_A{\ua^G},\bar\omega)&\cong \HH^1(G,\bar\omega)\cong \kk \\
\Ext^1_{\kG}(\kk_A{\ua^G},[\omega,\bar\omega])&\cong 
\HH^1(G,[\omega,\bar\omega])\cong \kk &
\Ext^1_{\kG}(\kk_A{\ua^G},[\bar\omega,\omega])&\cong \HH^1(G,[\bar\omega,\omega])\cong \kk 
\end{align*}
So if the number of summands isomorphic to any of these is infinite, then
$\Ext^1_{\kG}(\kk_A{\ua^G},M)$ woud be infinite dimensional.
The module $\kk_A{\ua^G}$ is uniserial $[\kk,\kk]$, and by the long 
exact sequence in $\Ext$, the degree one cohomology
$\HH^1(G,M)\cong\Ext^1_{\kG}(\kk,M)$ would be infinite dimensional,
contradicting our assumptions on $M$.
\end{proof}

\begin{lemma}\label{le:ll3}
The Loewy length of $M$ is at most three. Hence, $\Rad^2(M)\subseteq\Soc(M)$
and $\Rad(M)\subseteq \Soc^2(M)$.
\end{lemma}
\begin{proof}
The element $XYZ$ is in the socle of $\kG$, so it acts as zero on $M$
since $M$ is assumed to have no non-zero projective summands. The
image of $XYZ$ is the sum of the socles of the projective
indecomposables, and the quotient has Loewy length three.
\end{proof}

Let $d= d(M)$ be the total dimension of the non-projective part (the core) of
$M{\da_A}$. By Lemma~\ref{le:MdaA}, $d$ is finite. 
If $d(M)=0$ then by Proposition~\ref{pr:res-nocoho},
$M$ is a no-cohomology module and we are
done. So we assume by induction that $d(M)$ is positive and 
as small as possible.

\begin{lemma}\label{le:XYM}
We have $XYM=\Soc(M)$ and $\Ker(XY,M)=\Rad(M)$.
\end{lemma}

\begin{proof}
The inclusions $XYM\subseteq \Soc(M)$ and 
$\Rad(M)\subseteq \Ker(XY,M)$ follow  from
Lemma~\ref{le:ll3}. If there is an element of $\Soc(M)$ that is
not in $XYM$, then it generates a one dimensional submodule $M_1$
such that $d(M_1) = d(M)-1$ and we are done by induction.  
So we assume that $\Soc(M)
\subseteq XYM$. If there is an element of $\Ker(XY,M)$ 
not in $\Rad(M)$ then there is
a codimension one submodule $M_2$ of $M$ not containing it. 
Since $d(M_2) = d(M)-1$, we are again done by induction. 
So we have that $\Ker(XY,M)\subseteq \Rad(M)$.
\end{proof}

Now suppose that $U$ is a non-projective indecomposable summand $M{\da_A}$. 
Then $U\subseteq\Ker(XY,M)=\Rad(M)$ and so $XU\subseteq\Rad^2(M)\subseteq
\Soc(M)=XYM$. But $XYM\cap U=0$, so $XU=0$. Similarly $YU=0$, so $U$
is simple. Thus $U$ is not isomorphic to one of the uniserial modules 
$[\omega,\bar\omega]$ and $[\bar\omega,\omega]$ in the list \eqref{eq:A},
and $U$ must have dimension one. 

Suppose that the eigenvalue of $t$ on $U$ is $\omega$; the case where
the eigenvalue is $\bar\omega$ is similar. 
So $U$ is generated by an element $m$ with
$tm=\omega m$, $Xm=0$ and $Ym=0$.
Then $Zm$ is annihilated by $X$, $Y$ and $Z$, and
$tZm=\omega Zm$.
It follows that $Zm$  is in $\Soc(M)$, which equals $XYM$ by 
Lemma~\ref{le:XYM}.  Hence by Proposition~\ref{pr:eigenvalues},
there exists $m_1\in M$ such that $XYm_1=Zm$, and $tm_1=\omega m_1$.

Now $tZYm_1=\bar\omega ZYtm_1=ZYm_1$, and $ZYm_1$ is in
$\Soc(M)$. Since $\Hom_\kG(\kk,M)=0$, this implies that
$ZYm_1=0$.  Thus 
the linear span in $M$ of $m$, $Zm$ and $Ym_1$ is a $\kG$-submodule
of $M$ of dimension three.
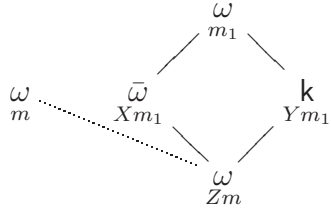
\[ \xymatrix@=3mm{
&&&\underset{m_1}\omega\ar@{-}[dl]\ar@{-}[dr]\\
\underset{m}\omega\ar@{.}[drrr]&&\underset{Xm_1}{\bar\omega}\ar@{-}[dr]
&&\underset{Ym_1}\kk\ar@{-}[dl]\\
&&&\underset{Zm}\omega} \]

We have $tXm_1=\omega Xtm_1=\bar\omega Xm_1$. Suppose that $Xm_1$ 
can be written as a linear combination $Ym_2+Zm_3$.
As in Proposition \ref{pr:eigenvalues} we may assume that $m_2$ and $m_3$
are eigenvectors of $t$ with eigenvalues $1$ and $\bar\omega$, respectively.
Then since $\Hom_\kG(M,\kk)=0$, we have $m_2\in \Rad(M)$. But then 
$Ym_2\in \Rad^2(M)\subseteq \Soc(M)$, so $ZYm_2=0$.
It follows that $ZXm_1=ZYm_2+Z^2m_3=0$.
Now let $M_1$ be the linear span of $m$, $Zm$, $Xm_1$ and $Ym_1$.
This is a $\kG$-submodule of $M$. Choose a submodule $M_2$ of $M$ of
codimension one, not containing $m_1$. Then $d(M_2/M_1) = d(M)-1$ and 
we are done by induction.

We are left with the case that $Xm_1$ cannot be written as a linear combination
$Ym_2+Zm_3$. In this situation, let $M_1$ be the linear span of 
$m$, $Zm$ and $Ym_1$, and $M_2$ be a codimension two submodule of 
$M$ not containing $m_1$ or
$Xm_1$. Then again $d(M_2/M_1) = d(M) -1$.
Hence, by induction, we are done. 

This completes the proof of Theorems~\ref{thm:tech} and~\ref{th:main}
for $G=\bZ/2\times A_4$. With the results of the last section, the proofs 
of the theorems are finished.

\bibliographystyle{amsplain}
\bibliography{../../repcoh}

\end{document}